\documentclass[a4paper,10pt,leqno]{amsart}
\usepackage[applemac]{inputenc}
\usepackage{epsfig,xypic}
\usepackage{amsthm,amssymb,bbm}%,a4wide}

\def\End{\mathrm{End}}
\newtheorem*{theo}{Theorem}

\def\k{\mathbbm{k}}

\def\R{\mathbbm{R}}

\author{Ivan Marin}
\title{On the residual nilpotence of pure Artin groups}
\date{January 2005}
\begin{document}

\maketitle

{\bf Abstract.} In this very short note we show that the residual
nilpotence of pure Artin groups of spherical type is easily
deduced from the faithfulness of the Krammer-Digne representations.

\bigskip

Let $W$ be a finite Coxeter group,
$B$ the associated (spherical) Artin group of spherical type, and $P$ the corresponding pure
group, defined as the kernel of the projection $\pi : B \twoheadrightarrow W$.
It is conjectured that the group $P$ is
residually nilpotent, meaning that the intersection of the terms
$C^r P$ for $r \geq 1$ of its descending central series, defined by
$C^1 P$ equals the commutator group $(P,P)$ and $C^{r+1} P = (P, C^r P)$,
is trivial. This fact is obvious for the dihedral types,
since it reduces in this case to the study of the free groups, which
is known for a long time. It
was proved by Falk and Randell for
types $A_n$, $B_n$ (and $I_2(m)$) in \cite{FALKDEUX} as an
easy consequence of their previous work \cite{FALK}, and by Markushevich
in \cite{MARK,MARKL} for type $D_n$.

There is a short way to prove this for every Artin group of 
spherical type using the generalization of the Krammer
representation by F. Digne (these representations where also obtained in type
ADE by A.M. Cohen and D.B. Wales in \cite{CW}) --- in the
same spirit that the linearity of these groups imply their
residual finiteness by Malcev theorem. It seems that this very simple
application
had been unnoticed yet.

\begin{theo} The pure Artin groups of spherical type are
residually nilpotent.
\end{theo}
\begin{proof}
First we may assume that $W$ is reduced. Moreover let's
assume first that $W$ is of crystallographic type. Let $N$ be the number of reflexions of $W$ and $V \simeq
\R^N$ the $\R$-vector
space formed by linear formal combination of reflexions of $W$. We
let $L = \R[q,q^{-1},t]$, $A = \k[[h]]$ where $q,t,h$ are indeterminates.
There are embeddings $\iota : L \hookrightarrow A$ given, for example, by $q \mapsto
e^h$ and $t \mapsto e^{\sqrt{2} h}$. We let $V_{q,t} = V \otimes_{\R} L$,
and $V_h = V\otimes_{\R} A = V_{q,t} \otimes_{\iota} A$. There
is a reduction morphism $\varphi : V_h \mapsto V$. F. Digne defined
in \cite{DIGNE} a faithful representation $R_W : B \to GL(V_{q,t}) \subset GL(V_h)$,
and it is easily checked from his formulas that the following
diagram is commutative
$$
\xymatrix{
B \ar[r] \ar[d] & GL(V_h) \ar[d]_{\varphi} \\
W \ar[r]_{R_W^0} & GL(V) }
$$
where the map $R_W^0 : W \to GL(V)$ is given by the natural conjugation action
of $W$ on its set of reflexions. It follows, using the identifications
$V = \R^N$, $\End(V_h) = M_N(A)$, that $R_W(C^r P) \subset 1 + h^{r+1}
M_N(A)$. In particular, if $x$ belongs to the intersection
of the terms $C^r P$, $R_W(x)$ belongs to the intersection
of the terms $1 + h^{r+1} M_N(A)$, which is trivial. Hence $R_W(x) = 1$
and $x = 1$ by faithfulness of the Krammer-Digne representation.

If $W$ is of type $H_3$, we let $W'$ be of type $D_6$, and denote
by $B'$, $P'$ the Artin group and pure Artin group associated to it.
There exists an injective morphism $j : B \to B'$ which sends $P$ to $P'$ and
induces an injective morphism $\check{j} : W \to W'$ such that the
following diagram commutes
$$
\xymatrix{
B \ar[r]^{j} \ar[d] & B'\ar[d] \\
W \ar[r]_{\check{j}} & W' }
$$
where the vertical arrows are the natural projections.
If $\sigma_1,
\sigma_2,\sigma_3$ denote the Artin generators of $B$ such that
$<\sigma_1,\sigma_2>$ is of type $I_2(5)$ and $<\sigma_2,\sigma_3>$
has type $I_2(3)$, and the Artin generators $\sigma'_1,\dots,\sigma'_6$
of $B'$ are numbered according to Bourbaki (see \cite{BOURB}), $j$ is defined
by $\sigma_1 \mapsto \sigma'_3 \sigma'_5$,
$\sigma_2 \mapsto \sigma'_2\sigma'_4$, $\sigma_3 \mapsto 
\sigma'_1\sigma'_6$. This is one of the ``folding morphisms'' defined
by J. Crisp in \cite{CRISP}. If $W$ is of type $H_4$, let $W'$ be of
type $E_8$, let $\sigma_1,
\sigma_2,\sigma_3, \sigma_4$ be the Artin generators of $H_4$ such
that $<\sigma_1,
\sigma_2,\sigma_3>$ is of type $H_3$ and number the Artin generators
$\sigma'_1,\dots,\sigma'_8$ of $B'$ according to Bourbaki. There is also
a folding morphism defined by
$\sigma_1 \mapsto \sigma'_3 \sigma'_5$,
$\sigma_2 \mapsto \sigma'_4\sigma'_6$, $\sigma_3 \mapsto 
\sigma'_3\sigma'_7$ and $\sigma_4 \mapsto \sigma'_1\sigma'_8$.
Then $j$ restricts to a monomorphism $P \to P'$, and $R_{W'}\circ j$
is faithful. Then $R_{W'}\circ j(P)
\subset 1 + h M_N(A)$, where $N$ is the number of reflexions of $W'$.
It follows that $P$ is residually nilpotent by the same argument.
Since there are also folding morphisms from the groups of type $I_2(m)$
to any Artin group of Coxeter number $m$, the conclusion also holds
for type $I_2(m)$ (alternatively one may use the faithfulness of the Iwahori-Hecke
representation proved in \cite{LEH}).

\end{proof}

Notice that we could have restrict ourselves to the representations
constructed by A.M. Cohen and D.B. Wales (see \cite{CW}) in type ADE since all Artin groups of
spherical type embed in one of these by a folding morphism.

\end{document}